\documentclass[letterpaper, 10 pt, conference]{ieeeconf}  

\IEEEoverridecommandlockouts                              

\usepackage{graphicx}
\usepackage{amsmath}
\usepackage{color}
\usepackage{epstopdf}
\usepackage{multirow}
\usepackage{amssymb}
\usepackage[table,xcdraw]{xcolor}
\usepackage{float}

\title{\LARGE \bf
Optimization-based Residential Load Scheduling to Improve Reliability in the Distribution Grid
}

\author{Abdulelah H. Habib, Elizabeth L Ratnam, Vahid R. Disfani, Jan Kleissl and Raymond A. de Callafon
\thanks{Abdulelah Habib, Elizabeth L Ratnam, Vahid R Disfani, Jan Kleissl and Raymond A. de Callafon are with the Department of Mechanical and Aerospace Engineering, University of California, San Diego
       {\tt\small {ahhabib,eratnam,disfani,jkleissl,callafon}
       @ucsd.edu}}
}
\begin{document}

\maketitle
\thispagestyle{empty}
\pagestyle{empty}

\begin{abstract}
\def\liz{\textcolor{blue}}
Despite the recent rapid adoption of rooftop solar PV for residential
customers, islanded operation during grid outages remains elusive for
most PV owners. In this paper we consider approaches to improve the reliability of electricity supply in the context of a residential microgrid, consisting of a group of residential customers each with rooftop solar PV, that are connected to the distribution network via a single point of common coupling. It is assumed that there is insufficient PV generation at all times to meet the electricity demand within the residential microgrid. Three optimization-based algorithms are proposed to improve the reliability of electricity supply to each residential customer, despite variability and intermittency of the solar resource and periods of infrequent and sustained power outages in the electricity grid. By means of a case study we show that the majority of residential customers achieve greater reliability of uninterrupted electricity supply when connecting to the residential microgrid in comparison to operating in isolated self-consumption mode.
\end{abstract}


\section{INTRODUCTION}
\def\paraph{\textcolor{blue}}
\def\vahid{\textcolor{red}}
\def\liz{\textcolor{blue}}
Large scale power outages or blackouts typically lead to millions of dollars in losses for industry, commercial and residential customers \cite{futureofmg,NYblackout,sandiegoblackout}. These power outages can be caused by human error, equipment failure, or may result from natural disasters such as the blackout induced by Hurricane Sandy in the Northeast of the U.S. in 2011 causing power shut off for 8 million customers for days and weeks with estimated damages around 50 Billion U.S. dollars \cite{sandyblackout,sandyblackstart}. 

A microgrid is designed to be interconnected with a medium voltage network under normal conditions, and to serve as a stable backup resource in case of isolation or islanding from the transmission grid (emergency operation mode)\cite{Coupledmicrogrids}. Forced isolation may occur in cases of voltage collapse, electric faults, or drops in power quality \cite{Coupledmicrogrids} at the point of common coupling (PCC). If the microgrid is well-designed, the transition to islanded operation should ideally occur smoothly with matching voltage and current phases on the PCC \cite{microgridemergency}.

Several authors have considered transition approaches to improve the power quality of the microgrid as it switches to island mode \cite{islandmode,hierarchicalcontrol,freqresponse}. 
Load scheduling is often a part of a transition approach to ensure power demand can be managed during islanding mode, given the limited power supply and bandwidth of the (renewable) distributed energy resources within the microgrid. Typically, large electric loads such as electric vehicle charging and household load appliances are considered in load scheduling, see e.g. \cite{EVscheduling} and \cite{Loadscheduling}. 

As load scheduling involves optimization to minimize electric energy losses (i.e. dumping of solar energy for lack of a load to utilize it) over a finite number of loads, numerical methods that exploit convex optimization routines in Mixed-Integer Linear Programming (MILP) are a promising method for optimal load scheduling. Applications of MILP in power systems can be seen in a variety of areas such as unit commitment of power production \cite{MILPpower}, power distribution network expansion \cite{MIPtran,MILPpowertrans}, scheduling of generation units in off-grid conditions in order to maximize supply performance of the system \cite{MNLPmicrogrid} as well as optimal scheduling of a renewable microgrid in an isolated load area \cite{Chung_microgrid}. MILP is also used in the optimal decentralized energy management problem of a microgrid \cite{MILPcdc}.



Based on the computational tools of MILP, this paper considers an optimization-based problem for scheduling loads of a group of residential customers, each with rooftop solar PV, that are connected to the distribution network via a single point of common coupling. Note that the problem formulation is identical for solar PV array on the roof of an apartment or condominium building, where the solar array could serve the load of the common areas as well as some of the units during a power outage. In the case of a single building the hardware (microgrid controller) and governance (building owner or homeowners association) issues are much more straightforward than for different buildings in a neighborhood. In the MILP formulation, residential participation (each house with a single meter is considered a load) is parametrized with a binary or integer component, while the optimization aims to maximize the number of loads that receive power, despite the insufficient PV power generation at all times to meet the electricity demand within the residential microgrid. The optimization is applied to residential grid data obtained from ten residential customers (houses) that were selected for analysis from the Australian grid from \cite{ratnam2015residential}. Different operating strategies for the members of the microgrid were considered including isolated self-consumption and several inter-connected sharing strategies to improve the reliability of supply to each residential customer in the microgrid during periods of infrequent and sustained power outages that result in isolation of the microgrid from the electric power network. 
The MILP optimization shows that a majority of the residential customers achieve greater reliability of electricity supply when connecting to the microgrid via "inter-connected sharing" in comparison to operating in isolated "self-consumption" mode.

Section~\ref{Preliminary} introduces the dataset used in this work and showcases microgrid simulation results for select houses for a high PV and low PV day. The problem formulation is discussed in Section \ref{problemformulaiton} with the mathematical optimization problem and the constraints. The simulation results that quantify the benefits of different operating strategies are covered in the Section \ref{results} and Conclusions follow.

\section{Problem Description}\label{Preliminary}
\subsection{Strategies for solar energy sharing in a microgrid}
\begin{figure}[ht]\centering
\includegraphics[width=.9\columnwidth]{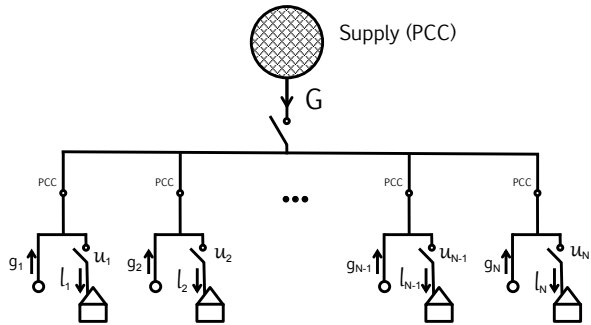}
\caption{Power generation and distribution configuration for a residential microgrid of ten houses with the possibility of islanding.}
\label{graphtheory}
\end{figure}

In this problem we assume a power system with bulk supply production as a generation unit connected to a distribution system via a main circuit breaker (CB) at the point of common coupling (PCC) to isolate the microgrid of residential customers each with their owns loads and PV system. Furthermore, each residential customer (also indicated exchangably by "house") has an additional CB referred to in Figure~\ref{graphtheory} as $u_i$, where $i$ is the house index. In cases of power outages (blackouts), faults or power quality disruption from the main power supply, the main CB at the PCC opens and, in addition, a certain set of houses decides to isolate itself to create a modified microgrid of residential customers. The decision making process that decides which houses to (dis)connect is managed by an optimization problem that operates at a time step of 30 minutes.

We select a subset of a dataset with 300 de-identified residential customers with PV in a distribution network in Australia. The optimization problem considers one year to cover many load and PV scenarios that may occur within a year for ten houses. The PV systems vary in size. The total PV rated capacity is 17.33 kW$_{\rm AC}$. Daily peak solar power averages around 11 kW. The corresponding daytime load peak is around 8 kW with a higher peak in the evening that reaches 13 kW. For model validation, the first ten houses were selected with customers ID [2 13 14 20 33 35 38 39 56 69] for July 1, 2010 through June 30, 2011 defined in \cite{ratnam2015residential}.

The two main microgrid operational modes considered in this case study are "isolated self-consumption" and "inter-connected sharing mode". 
Isolated self-consumption is the case where each house has been disconnected and can only be supply its loads from its own solar power. If the load is higher than the solar power at any time step $t_k$, then no power is supplied and the solar power is considered to be lost.
 In inter-connected sharing mode, houses will exchange PV power to supply their electric loads. Within the inter-connected sharing mode, three different sub strategies are investigated:

\begin{itemize}

\item Strategy A is to maximize the number of customers to be supplied.
\item Strategy B is to maximize the time duration of supplied load or what will be referred to as number of switches.
\item Strategy C is to minimize losses due to unutilized solar energy.

\end{itemize}

\noindent
Additional constraints are added for all strategies and include a minimum up-time and down-time for a supplied load event. Additional conditions will be included such as a "fairness weighting matrix" where customers are prioritized based on certain criteria. In this case study we use the percentage of PV self-generation with respect to the load of each house as a "fairness weighting".

\subsection{Illustration of Data}
\begin{figure}[ht]\centering
\includegraphics[width=.9\columnwidth]{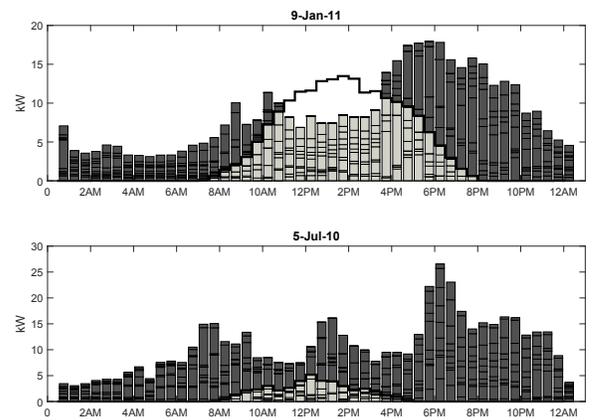}
\caption{Two sample days from a day in summer (top) and winter (bottom) where the dark gray bars represent total load and the light gray bars shows load that could be supplied from solar energy in case of an islanded microgrid. 
Each sub bar shows a different house and the total shows the aggregation of bar stacks. Note that January is summer in the southern hemisphere. All results later consider 1 year.}
\label{samplesdays}\vspace{-10pt}
\end{figure}

Two sample days illustrate extremes in potential for solar energy to power the microgrid (Figure~\ref{samplesdays}). These results are presented here to guide the problem formulation. 
In summer, solar generation is high compared to load and on this specific day it happens to exceeds load at the solar peak. On the other hand, higher loads in winter correlate with low solar generation. The load supplied by solar is what solar was able to supply for each house in the isolation mode where the total load was not met by solar for the whole system. 
\begin{figure}[ht]\centering
\includegraphics[width=1.0\columnwidth]{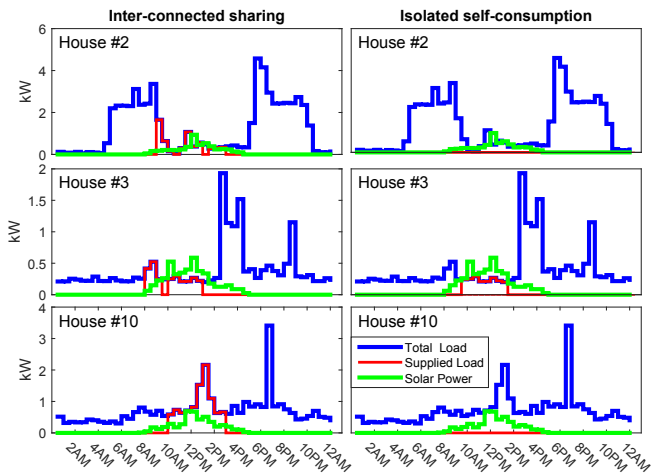}
\caption{Hourly results for three different customers for the winter day, where the blue line shows total load, the red line represents supplied load, and the green line represents solar power for each house. The results are shown for inter-connected sharing and isolated self-consumption case.
}
\label{samplehouses}\vspace{-15pt}
\end{figure}

Figure~\ref{samplehouses} shows the consumption from three select houses and solar generation behavior for inter-connected sharing and isolated self-consumption operational mode. For this case, the minimum up-time and down-time are 1.5 h (3 time steps). 
The isolated self-consumption operational mode does not allow solar generation to supply any load unless the solar generation exceeds the load for the minimum up-time constraint. Furthermore, if the house is "turned off" (disconnected), it has to be off for minimum down-time. The minimum up-time effect occurs in some houses in Figure~\ref{samplehouses}: for example house \# 2 around noon and around 2PM where PV generation was temporarily higher than load, but no load was supplied to the house due to the minimum up-time constraint. Self-consumption is therefore only attractive if the size of the PV associated with each house is large relative to the consumption, but most houses, at least in winter, do not get any power from their local PV generation and all solar generation is lost. 
On the other hand, the inter-connected sharing operation mode can aggregate all PV generation to serve more customers. Allocations of energy to particular houses must be determined based on an optimization to maximize solar energy utilization and/or customer supply.  Figure~\ref{samplehouses} shows that some houses, at least during parts of the day, enjoyed load supply even though load exceeded PV self-generation. More quantitative results will be shown and discussed in Section~\ref{results}.

\section{Problem Formulation}
\label{problemformulaiton}
The optimization problem for "inter-connected sharing mode" to address the three different sub strategies A, B and C defined in Section~\ref{Preliminary} considered in this paper can be presented in the general format of
\begin{equation}
\begin{aligned}
& \max_\mathbf{U} &&\mathbf{f(U)}\\
& \text{s.t.}&&\mathbf{g(U)}\le0\\
&&& \mathbf{h(U)}=0\\
&&& \mathbf{U}\subset\{0,1\}^{|\mathbf{U}|}
\end{aligned}
\label{mainopt}
\end{equation}
where $\mathbf{f}$, $\mathbf{g}$ and $\mathbf{h}$ denote the mathematical formulations of the objective function, inequality constraints, and equality constraints respectively. The binary decision variable $\mathbf{U}=[u_1 ,u_2,\cdots,u_N]^T$ is the matrix of switching statuses for all houses at all time steps during a day such that $u_i$ is the column vector of switching statuses for house $i$ during one specific day.

\subsection{Objective Functions}
Two different objective functions are considered in this paper. The following equations describe $\mathbf{f_1(U)}$ which minimizes PV energy loss (maximizes load supply) and $\mathbf{f'_2(U)}$ which maximizes the number houses that are switched on, respectively: 
\begin{equation}
\begin{aligned}
&\mathbf{f_1(U)}=\mathbf{1}^T\mathbf{\cdot U\cdot1}
\end{aligned}
\label{f1}
\end{equation}
\begin{equation}
\begin{aligned}
&\mathbf{f'_2(U)}=\mathbf{1}^T\mathbf{\cdot(U\circ L)\cdot1}
\end{aligned}
\label{f2}
\end{equation}
where $\mathbf{1}$ is a column vector with the appropriate size, whose elements are all equal to one. 
The load matrix is also denoted by $\mathbf{L}=[l_1 ,l_2,\cdots,l_N]^T$ where $l_i$ is the column load vector of house $i$ during a day. The notation $\mathbf{U\circ L}$ is used to show the Hadamard (element-wise) product of the two matrices $\mathbf{L}$ and $\mathbf{U}$.
By defining $\mathbf{Y}=\mathbf{U}\circ \mathbf{L}$ as a new constraint, the objective function $\mathbf{f'_2}$ can be updated as,
\begin{equation}
\begin{aligned}
\mathbf{f_2(U)}=\mathbf{1}^T\mathbf{(Y)\cdot1}
\end{aligned}
\label{f3}
\end{equation}
The first objective function $\textbf{f}_1$ considers the number of switches which determines how many houses are supplied with power. Multiplying $\textbf{U}$ by a column vector of ones from both sides is another way of representing the sum over all houses of the sum of switching events for each house. 
This ensues that the number of houses who receive power at some point is maximized follows a utilitarian philosophy, but in doing so $\textbf{f}_1$ does not maximize the solar energy utilization. $\textbf{f'}_2$ aims to increase the energy supplied for the whole microgrid and it is represented as in $\textbf{f}_2$ by introducing a new variable $\textbf{Y}$.
\subsection{Constraints}
\subsubsection{Definition of $\mathbf{Y}$}
The substitution $\mathbf{U\circ L}$ by $\mathbf{Y}$ in the objective function motivates the following constraint of the supplied load matrix:
\begin{equation}
\begin{aligned}
\mathbf{Y}-\mathbf{U\circ L}=0
\end{aligned}
\label{const_Y}
\end{equation}

\subsubsection{Available Power}
To prevent frequency issues, the maximum total load that the microgrid can supply must be less than the total PV energy available at each time interval:
\begin{equation}
\begin{aligned}
\mathbf{1}^T\mathbf{\cdot(U\circ L)}\le\mathbf{1}^T\mathbf{\cdot G},
\end{aligned}
\label{const_power}
\end{equation}
where $\mathbf{G}=[g_1,g_2,\cdots,g_N]^T$ is the PV generation matrix. 

\subsubsection{Minimum Up-time and Minimum Down-time} To avoid damage to load units and inconvenience to residents because of frequent start-ups and shut-downs, a set of constraints are defined to guarantee that the unit is switched on (off) for a at least $m^+$ ($m^-$) time steps before it is switched off (on). These constraints are called minimum up (down) time and are defined as:
\begin{equation}
\begin{aligned}
u_{i,t_k}-\sum_{h=t_k-m^+_i +1}^{t_k} v_{i,h}\leq 0 && \forall_{m^+_i\le t_k\le T}\\
(1-u_{i,t_k})-\sum_{h=t_k-m^-_i +1}^{t_k} w_{i,h}\leq 0&& \forall_{m^-_i\le t_k\le T},
\end{aligned}
\label{const_minupdown}
\end{equation}
where the matrix $\mathbf{V}\subset\{0,1\}^\mathbf{|V|}$ and $\mathbf{W}\subset\{0,1\}^\mathbf{|W|}$ are denoted as start-up and shut-down matrices respectively, and their elements are defined as:
\begin{equation}
\begin{aligned}
&v_{i,t_k}-w_{i,t_k}=u_{i,t_k}-u_{i,t_{k-1}} &&\forall_{1\le i\le N}\forall_{2\le t_k\le T}\\
&v_{i,t_k}+w_{i,t_k}\leq 1 && \forall_{1\le i\le N}\forall_{2\le t_k\le T}\\
&v_{t,1}=w_{i,1}=0&& \forall_{1\le i\le N}
\end{aligned}
\label{const_vw}
\end{equation}

\subsubsection{Minimum Daily Connection to Grid}
The following constraint guarantees that each house is connected to the grid at least one minimum up-time. 
\begin{equation}
\begin{aligned}
\mathbf{1- U\cdot 1\le 0}
\end{aligned}
\label{const_minconnect}
\end{equation}

All the possible objective functions and constraints in \eqref{mainopt} through \eqref{const_minconnect} are linear. Thus the optimization problem is convex and can be solved via Mixed-integer Linear Programming (MILP) tools such as Gurobi using CVX. There exist many mature MILP solvers which are capable of solving large-scale MILP problems with millions of variables within a reasonable time frame \cite{MIPtalg}.

\subsection{Operational Strategies and Reasoning Indices}
\label{OperationStrategies}
There are three main operational strategies for inter-connected sharing and each strategy is associated with two sub-strategies. In strategy A, the goal is to maximize the objective function $\mathbf{f_1}$ with all constraints \eqref{const_Y}-\eqref{const_minconnect}, which maximizes the number of households whose load is served at some point. This forces all houses to receive power for at least one minimum up-time.
Strategy B also maximizes the objective function $\mathbf{f_1}$ while the minimum daily connection constraint in \eqref{const_minconnect} is neglected. The goal is to increase the number of switches without enforcing that all houses receive power at least once. 
Strategy C maximizes solar energy utilization by maximizing the objective function $\mathbf{f_2}$ considering all constraints other than \eqref{const_minconnect}. Solar energy is distributed in every possible way to reduce any losses even if it means that more houses never receive power.

The sub operation strategies A+, B+ and C+ are exactly the same as strategies A, B and C respectively with the following modifications to the objective functions:
\begin{equation}
\begin{aligned}
&\mathbf{f_1}^*=\mathbf{1}^T\mathbf{\cdot\mathbb{W}\circ U\cdot1}
\end{aligned}
\label{f1*}
\end{equation}
\begin{equation}
\begin{aligned}
\mathbf{f_2}^*=\mathbf{1}^T\mathbf{\cdot \mathbb{W}\circ Y\cdot1}
\end{aligned}
\label{f3*}
\end{equation}
where $\mathbb{W}$ is a "fairness weighting" matrix of the same size of $\textbf{U}$ and $\textbf{Y}$. The weighting matrix introduces preferential weighting for certain houses to receive power even though it deviates from the solutions for $\mathbf{f_1}$ and $\mathbf{f_2}$. The weights would be set by a governing entity based on perceived fairness criteria such as prioritizing houses with larger PV generation, lower load demand, or prioritizing critical loads (e.g. medical needs) either permanently or during a given time span. The weights could also be based on a market where individual homes pay to receive priority for load. For illustrative purposes, the ratio of PV generation divided by the total load is used as weighting function in $\mathbb{W}$ here. 

The six strategies are compared to the self-consumption strategy. This strategy can be modeled by solving the same optimization problem as in strategy C for each house where the number of houses in the problem is equal to $1$.

To compare the simulation results of different objective functions with the isolated self-consumption operational mode, two indices are defined in this paper. The first index is the percentage of supplied load, defined as below,
\[
\% \text{ of Load met}=(\textbf{Y}\cdot \textbf{1})\oslash (\textbf{L}\cdot \textbf{1})
\]
where $\oslash$ is defined as element-wise division.

Since $\textbf{Y, L}$ and $\textbf{G}$ are of size (N, T) multiplying  $\textbf{Y}$ by vector $\textbf{1}$ of size (T,1) results in (N,1) which is the sum of supplied load for each house for a specific day. The percentage of load met is the ratio of supplied load $\textbf{Y}$ divided by the total load $\textbf{L}$ for a specific day.

The other index is the percentage of PV utilization of a given day, which is determined as follows,
\[
\% \text{ PV Generation Utilization}=(\textbf{1}^T\cdot\textbf{G}) \oslash (\textbf{1}^T\cdot\textbf{L})
\]
and reports what portion of individual houses' solar generation is utilized by each house for the isolated self-consumption operational mode and what portion of total solar generation in microgrid is utilized in for the inter-connected sharing strategies. 

\section{Simulation Results }
\label{results}

\begin{figure}[ht]\centering
\includegraphics[width=1.0\columnwidth]{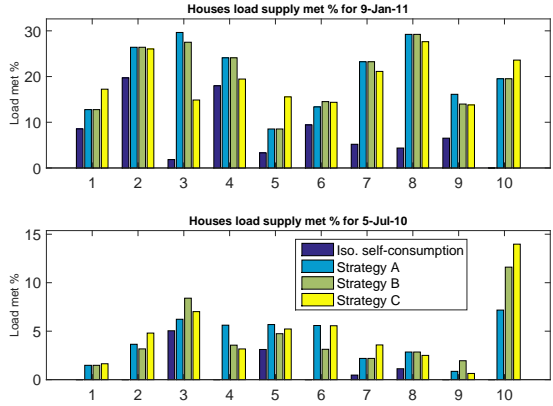}
\caption{Aggregated daily results for all 10 houses comparing different objective functions and isolated self-consumption based on percentage of load met for the summer and winter day of Figure 2.}
\label{onedayresults}
\end{figure}

Figure \ref{onedayresults} shows aggregate results for all 10 houses over the course of two sample days, based on self-consumption with the three sharing strategies A, B and C defined in Section \ref{OperationStrategies}. On the summer day around 20\% of loads are met for all house and all house are powered at least once during the day. The January 9 results are representative for summer days when the peak of the aggregate solar generation is higher than the total load and therefore all houses can be powered during that time, independent of strategy. Some houses did not benefit as much as others from energy sharing, for example house \# 3, 8, and 10 benefit the most while houses \# 4 or 6 only receive marginally more energy compared to self-consumption. House \# 10 was not able to self-consume any time during the day as the load was always higher than the solar generation. 

Comparing the inter-connected sharing operational modes A and B resulted in identical results for 6 out of 10 houses. In 3 houses strategy C met less load compared to A and B, while 2 houses had higher load met. Little to no change is expected since constraint (9) is satisfied for both strategies through the solar energy excess at midday. Strategy C results were also mostly similar to A and B, but only three houses (\# 1, 5, 10) benefited, while two houses (\# 3, 4) lost a substantial amount of energy.

July 5 is representative of a winter day with low solar power and high load demand which tends to emphasize differences between the strategies. Gains from sharing compared to self-consuming were larger: all 10 houses were powered with energy sharing while only 4 houses received some energy from their own solar generation. Results for the winter day also vary by operational strategy A, B, and C. It is interesting that only four houses (\# 1, 2, 7, 10) benefited from strategy C. Considering both days, only house \# 10 benefited from strategy C consistently. For these two days, there is no clear winner between the inter-connected sharing operational strategy, but it is clear that sharing energy is advantageous for every house. In the following, the results were analyzed for different months and one year to quantify the performance of each strategy based on certain objective functions.


\begin{figure}[ht]\centering
\includegraphics[width=.9\columnwidth]{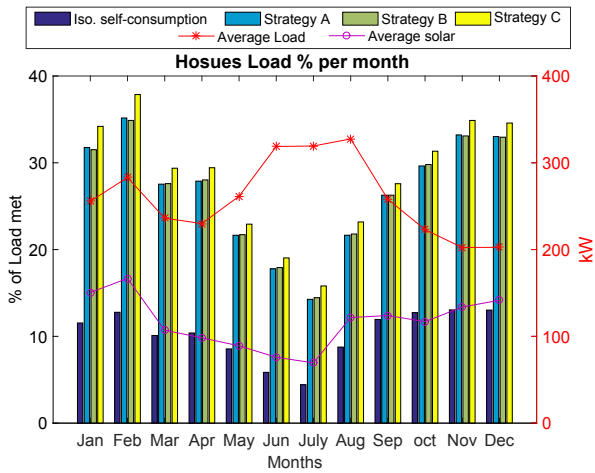}
\caption{Monthly results for percentage of load met for different strategies averaged over all houses. On the right hand axis, average load and solar generation for each month is presented.}
\label{monthly}
\end{figure}

Figure \ref{monthly} summarizes aggregate monthly average results for isolated self-consumption and the main three inter-connected energy sharing strategies. The seasonal (summer and winter) trends are consistent with results for the sample days. During summer time (October through February) there is a higher percentage of load met due to larger solar generation as well as load behavior. In all months, isolated self-consumption scores the lowest load met percentage while strategy C scores the highest with negligible differences between strategies A and B. 

\begin{figure}[ht]\centering
\includegraphics[width=.9\columnwidth]{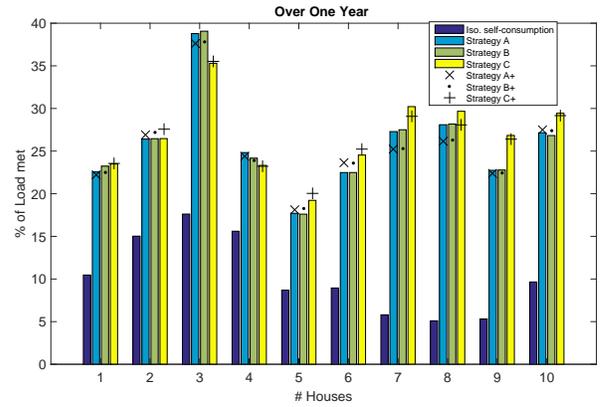}
\caption{Annual results for different operating strategies and houses.}
\label{housesyearly}\vspace{-15pt}
\end{figure}

Strategies A+, B+, and C+ use the same objective function, but include the weighting $\mathbb{W}$ that prioritizes additional constraints. The weighting matrix in this case is defined based on the ratio of solar generation to the load. Figure \ref{housesyearly} 
presents the annual load met percentage of each house for all six scenarios and the isolated self-consumption. All houses benefit from inter-connected energy sharing. This value varies over the houses but it is consistent that all are better in inter-connected sharing mode. 

Comparing different strategies, seven houses have higher load met percentage in strategy C, two houses have higher load met percentage in strategies A and B, and only one house with all three strategies equal. 
 The weighting strategies follow the pattern of the original ones for all houses. House \# 6 and 7 stand out, where in house \# 6 C+ strategy resulted in higher average load met percentage over a year opposite to house \# 7. 
Three houses benefit from applying the fairness weighting matrix while three are not affected, and the four others suffered. Among the three main strategies, C and C+ combined for the highest percentage of load met with seven out of ten houses (four and three houses respectively). 
\begin{table}[ht]
\scriptsize
\centering
\caption{\# of Houses supplied for different strategies.}
\label{numbofhouses}
\begin{tabular}{|c|c|c||c|c|c||c|c|c|}
\hline
\multicolumn{3}{|c||}{Strategy A}                                                                             & \multicolumn{3}{c||}{Strategy B}                                    & \multicolumn{3}{c|}{Strategy C}                                                                              \\ \hline
H                        & Days                     & \begin{tabular}[c]{@{}c@{}}\% Load \\ met\end{tabular} & H  & Days & \begin{tabular}[c]{@{}c@{}}\% Load \\ met\end{tabular} & H                        & Days                     & \begin{tabular}[c]{@{}c@{}}\% Load \\ met\end{tabular} \\ \hline
10                       & 359                      & 25.66                                                  & 10 & 341  & 26.07                                                  & 10                       & 347                      & 27.47                                                  \\ \hline
\cellcolor[HTML]{C0C0C0} & \cellcolor[HTML]{C0C0C0} & \cellcolor[HTML]{C0C0C0}                               & 9  & 352  & 25.55                                                  & 9                        & 356                      & 26.91                                                  \\ \hline
\cellcolor[HTML]{C0C0C0} & \cellcolor[HTML]{C0C0C0} & \cellcolor[HTML]{C0C0C0}                               & 8  & 359  & 25.07                                                  & 8                        & 363                      & 26.43                                                  \\ \hline
\cellcolor[HTML]{C0C0C0} & \cellcolor[HTML]{C0C0C0} & \cellcolor[HTML]{C0C0C0}                               & 7  & 362  & 24.83                                                  & 6                        & 365                      & 26.30                                                  \\ \hline
\cellcolor[HTML]{C0C0C0} & \cellcolor[HTML]{C0C0C0} & \cellcolor[HTML]{C0C0C0}                               & 6  & 365  & 24.66                                                  & \cellcolor[HTML]{C0C0C0} & \cellcolor[HTML]{C0C0C0} & \cellcolor[HTML]{C0C0C0}                               \\ \hline
\end{tabular}\vspace{-8pt}
\end{table}

Now we consider during how many days of the year a given number of houses receive at least power for the minimum  up-time for different operating strategies (Table \ref{numbofhouses}). As expected strategy A yields the most (359) days for all houses to be on; only for six days the optimization did not yield any results, i.e. the load was always larger than PV generation. The average percentage of load met for one year is 25.66\%. Relaxing the up-time constraint to 1 hour instead of (1.5 h and 0.5 h) allows all houses to be on for every day of the year. Strategy C was able to keep ten houses on for 347 days and six houses were on for the whole year. Strategy B scored the lowest in that all ten houses were on only for 341 days. The highest percentages of load met was achieved by strategy C.



\begin{table}[ht]
\centering
\scriptsize
\caption{Percent load met, percent PV generation utilized, and average number of houses supplied at least once for different operating strategies averaged over the year.}
\label{objfuncomp}
\begin{tabular}{c|c|cccccc|}
\cline{2-8}
\multirow{2}{*}{}                                                                            & \multirow{2}{*}{\begin{tabular}[c]{@{}c@{}}Iso. self \\ cons.\end{tabular}} & \multicolumn{6}{c|}{Strategies}                                                                                                                        \\
                                                                                             &                                                                             & A                          & A+                         & B                          & B+                         & C                          & C+    \\ \hline                          \hline                                           \multicolumn{1}{|c||}{\begin{tabular}[c]{@{}c@{}}\% of Load\\  met\end{tabular}}    & 9.17                                                                        & \multicolumn{1}{c|}{24.64} & \multicolumn{1}{c|}{24.43} & \multicolumn{1}{c|}{24.66} & \multicolumn{1}{c|}{24.49} & \multicolumn{1}{c|}{26.30} & 26.26 \\ \hline
\multicolumn{1}{|c||}{\begin{tabular}[c]{@{}c@{}}\% PV gen.\\  utilized\end{tabular}}                                                         & 29.25                                                                       & \multicolumn{1}{c|}{78.59} & \multicolumn{1}{c|}{77.93} & \multicolumn{1}{c|}{78.66} & \multicolumn{1}{c|}{78.11} & \multicolumn{1}{c|}{83.90} & 83.77 \\ \hline
\multicolumn{1}{|c||}{\begin{tabular}[c]{@{}c@{}}Avg \# \\houses to \\be supplied\end{tabular}} & 7.78                                                                        & \multicolumn{1}{c|}{10*}   & \multicolumn{1}{c|}{10*}   & \multicolumn{1}{c|}{9.70}  & \multicolumn{1}{c|}{9.69}  & \multicolumn{1}{c|}{9.61}  & 9.71  \\ \hline
\end{tabular}\vspace{-7pt}
\end{table}

Table \ref{objfuncomp} summarizes the performance of both operational mode with all strategies for one year. Isolated self-consumption is the worst operational mode in terms of load met percentages and PV utilization where it scored 65\% less than the best strategy. In terms of average number of houses to be supplied isolated self-consumption also scored the lowest. C+ and C strategies differ by less than 0.2\% and score the highest percentage of load met and PV utilization. In general, a strategy and its weighted version (for example C and C+) are expected to yield similar results, since the weighted strategy only changes the priority of which house is supplied. 
Moreover, the computational time for the optimization using MATLAB and CVX (Gurobi solver) was less than a second performed in a 3.4 GHz Intel Core i7 processor with 32 GB of RAM. 

\section{Conclusions}
In this paper we propose optimization-based residential customer scheduling to improve the reliability of electricity supply for residential customers during islanded microgrid operation. Each residential customer owns rooftop PV that can be used to supply just their own load or be shared across the microgrid to satisfy different operational strategies. In the latter case, residential customer scheduling is based on mixed-integer linear programming in which integers are used to parametrize the power status of each house and linear constraints enforce minimum up-time and down-time of power provision. The different operating strategies to distribute PV energy across the members of the microgrid include different objective functions which focus on the optimal use of solar PV within a microgrid: A) Forcing all houses to receive power at least once, B) Maximizing the number of switches without forcing all houses to be connected, C) maximizing the utilization of available solar power distributed among the grid to reduce power losses. Additional strategies were considered which used a priority or fairness weighting matrix to determine scheduling. The weighting matrix was computed by considering the load-to-generation ratio for each house, but other weighting based on priority of the loads in each house can be considered. 

A case study based on historical yearly data for ten houses was conducted. The mixed-integer linear programming results show that isolated self-consumption operation was the worst option for all houses. The objective which maximized the use of available solar power resulted in the highest percentage load met. Although results vary for each house, the trends over the year are consistent. Future work will include backup generation such as storage and distributed energy resources.
\bibliographystyle{IEEEtran}

\end{document}